\RequirePackage{ifpdf}
\ifpdf 
\documentclass[pdftex]{sigma}
\else
\documentclass{sigma}
\fi

\newcommand \al{\alpha}
\newcommand\be{\beta}
\newcommand\ga{\gamma}
\newcommand\de{\delta}
\newcommand\ep{\varepsilon}

\newcommand\et{\eta}
\renewcommand\th{\theta}

\newcommand\ka{\kappa}

\newcommand\rh{\rho}

\newcommand\si{\sigma}

\newcommand\ph{\varphi}

\newcommand\ps{\psi}
\newcommand\om{\omega}
\newcommand\Ga{\Gamma}
\newcommand\De{\Delta}

\newcommand\La{\Lambda}

\newcommand\Om{\Omega}

\newcommand\resp{resp.\ }
\newcommand\ie{i.e.\ }
\newcommand\wrt{w.r.t.\ }

\newcommand\X{\mathfrak X}
\newcommand\on{\operatorname}
\newcommand\g{\mathfrak g}
\newcommand\h{\mathfrak h}
\newcommand\kkk{\mathfrak k}
\newcommand\gl{\mathfrak gl}

\renewcommand\div{\on{div}}
\newcommand\ad{\on{ad}}

\newcommand\tr{\on{tr}}
\newcommand\Ad{\on{Ad}}
\newcommand\grad{\on{grad}}
\newcommand\Flux{\on{Flux}}
\newcommand\flux{\on{f\/lux}}
\def\RR{\mathbb R}
\newcommand\Diff{\on{Dif\/f}}
\newcommand\SDiff{\on{SDif\/f}}
\newcommand\Ric{\on{Ric}}
\newcommand\Der{\on{Der}}
\newcommand\Def{\on{Def}}
\newcommand\x{\times}
\newcommand\oo{\infty}
\newcommand\ZZ{\mathbb Z}
\newcommand\TT{\mathbb T}
\newcommand\curl{\on{curl}}

\begin{document}
\numberwithin{equation}{section}


\renewcommand{\PaperNumber}{030}

\FirstPageHeading

\renewcommand{\thefootnote}{$\star$}

\ShortArticleName{Geodesic Equations on Dif\/feomorphism Groups}

\ArticleName{Geodesic Equations on Dif\/feomorphism Groups\footnote{This paper is a contribution to the Proceedings
of the Seventh International Conference ``Symmetry in Nonlinear
Mathematical Physics'' (June 24--30, 2007, Kyiv, Ukraine). The
full collection is available at
\href{http://www.emis.de/journals/SIGMA/symmetry2007.html}{http://www.emis.de/journals/SIGMA/symmetry2007.html}}}

\Author{Cornelia VIZMAN}

\AuthorNameForHeading{C. Vizman}

\Address{Department of Mathematics, West University of Timi\c soara, Romania}
\Email{\href{mailto:vizman@math.uvt.ro}{vizman@math.uvt.ro}}

\ArticleDates{Received November 13, 2007, in f\/inal form March
01, 2008; Published online March 11, 2008}

\Abstract{We bring together those systems of hydrodynamical type
that can be written as geodesic equations on dif\/feomorphism groups
or on extensions of dif\/feomorphism groups with right invariant $L^2$ or $H^1$ metrics. We present their formal derivation starting from Euler's equation, the f\/irst order equation satisf\/ied by the right logarithmic derivative of a~geodesic in Lie groups with right invariant metrics.}


\Classification{58D05; 35Q35}


\begin{flushright}\it
A fluid moves to get out of its own way
as efficiently as possible.\\[1mm]
Joe Monaghan
\end{flushright}


\section{Introduction}

Some conservative systems of hydrodynamical type
can be written as geodesic equations on
the group of dif\/feomorphisms or the group of volume preserving dif\/feomorphisms of a
Riemannian
manifold, as well as on extensions of these groups.
Considering right invariant $L^2$ or $H^1$ metrics on these inf\/inite dimensional Lie
groups,
the following geodesic equations can be obtained:
the Euler equation of motion of a perfect f\/luid \cite{Arnold, EM},
the averaged Euler equation \cite{MRS,Shkoller},
the equations of ideal magneto-hydrodynamics \cite{VD,MRW},
the Burgers inviscid equation \cite{Burgers}, the template matching equation
\cite{HMA,V3},
the Korteweg--de Vries equation \cite{OK},
the Camassa--Holm shallow water
equation \cite{CH,Misiolek1,Kouranbaeva},
the higher dimensional Camassa--Holm equation
(also called EPDif\/f or averaged template matching equation)
\cite{HM}, the superconductivity equation \cite{Roger},
the equations of motion of a charged ideal f\/luid \cite{V1},
of an ideal f\/luid in Yang--Mills f\/ield \cite{GR3}
and of a stratif\/ied f\/luid in Boussinesq approximation \cite{Zeitlin2,V2}.

For a Lie group $G$ with right invariant metric, the geodesic equation written for
the right logarithmic derivative $u$ of the geodesic is a f\/irst order equation
on the Lie algebra $\g$, called the {\it Euler equation}. Denoting by $\ad(u)^\top$
the adjoint of $\ad(u)$ with respect to the scalar product on $\g$ given by the metric,
Euler's equation can be written as $\frac{d}{dt}u=-\ad(u)^\top u$.
In this survey type article we do the formal derivation of all the equations of
hydrodynamical type
mentioned above, starting from this equation.

By writing such partial dif\/ferential
equations as geodesic equations on dif\/feomorphism groups,
there are various properties one can obtain using the Riemannian geometry of right
invariant metrics on these dif\/feomorphism groups.
We will not focus on them in this paper, but we list some of them below, with some
of the references.

For some of these equations smoothness of the geodesic spray on the group implies local
well-posedness
of the Cauchy problem as well as smooth dependence on the initial data. This applies for
the following right invariant Riemannian metrics: $L^2$ metric on the group of volume
preserving
dif\/feomorphisms \cite{EM}, $H^1$ metric on the group of volume preserving dif\/feomorphisms
on a boundary free manifold \cite{Shkoller},
on a manifold with Dirichlet boundary conditions \cite{MRS,Shkoller2} and with Neumann or mixt
boundary conditions \cite{Shkoller2,GR}, $H^1$ metric on the group of dif\/feomorphisms of the
circle \cite{Shkoller,Kouranbaeva} and on the Bott--Virasoro group \cite{CKKT},
and $H^1$ metric on the group of dif\/feomorphisms on a higher dimensional manifold
\cite{GR2}.

There are also results on the sectional curvature (with information  on the Lagrangian
stability)
\cite{Arnold,NHK,Preston,Misiolek2,MR,PS,V6,Hattori,ZK,ZP,V1},
on the existence of conjugate points \cite{Misiolek,Misiolek3}
and minimal geodesics \cite{Brenier}, on the f\/initeness of the diameter
\cite{Shnirelman,Sh,ER}, on the vanishing of geodesic distance \cite{MM},
as well as on
the Riemannian geometry of subgroups of dif\/feomorphisms as a submanifold of the full
dif\/feomorphism group \cite{Misiolek0,BR,KM3,V3}.


\section{Euler's equation}

Given a regular Fr\'echet--Lie group in the sense of Kriegl--Michor
\cite{KM}, and a (positive def\/inite) scalar product
$\langle \ , \ \rangle :\g\times\g\to\RR$ on the Lie algebra $\g$, we can def\/ine a right
invariant metric on $G$ by
$g_x(\xi,\eta)=\langle\xi x^{-1},\eta x^{-1}\rangle$
for $\xi,\eta\in T_xG$.
The energy functional of a smooth curve $c:I=[a,b]\to G$ is def\/ined by
\begin{gather*}
E(c) =\frac12\int_a^b g_{c(t)}(c'(t),c'(t))dt
=\frac12\int_a^b\langle \de^rc(t),\de^rc(t)\rangle dt,
\end{gather*}
where $\delta^r$ denotes the right logarithmic derivative (angular velocity) on the
Lie group $G$, \ie $\de^rc(t)=c'(t)c(t)^{-1}\in\g$.
We assume the adjoint of $\ad(X)$ with respect to $\langle \ , \ \rangle $ exists for all $X\in\g$
and we denote it by $\ad(X)^\top$, \ie
\begin{gather*}
\langle \ad(X)^\top Y,Z\rangle =\langle Y,[X,Z]\rangle, \qquad\forall \, X,Y,Z\in\g.
\end{gather*}
The corresponding notation in \cite{AK} is $B(X,Y)=\ad(Y)^\top X$
for the bilinear map $B:\g\x\g\to\g$.


\begin{theorem}\label{theo}
The curve $c:[a,b]\to G$ is a geodesic for the right invariant metric $g$ on $G$
if and only if its right logarithmic derivative $u=\de^rc:[a,b]\to\g$
satisfies the Euler equation:
\begin{gather}\label{euler}
\frac{d}{dt}u=-\ad(u)^\top u.
\end{gather}
\end{theorem}
\begin{figure}[h]
\centerline{\includegraphics[width=97mm]{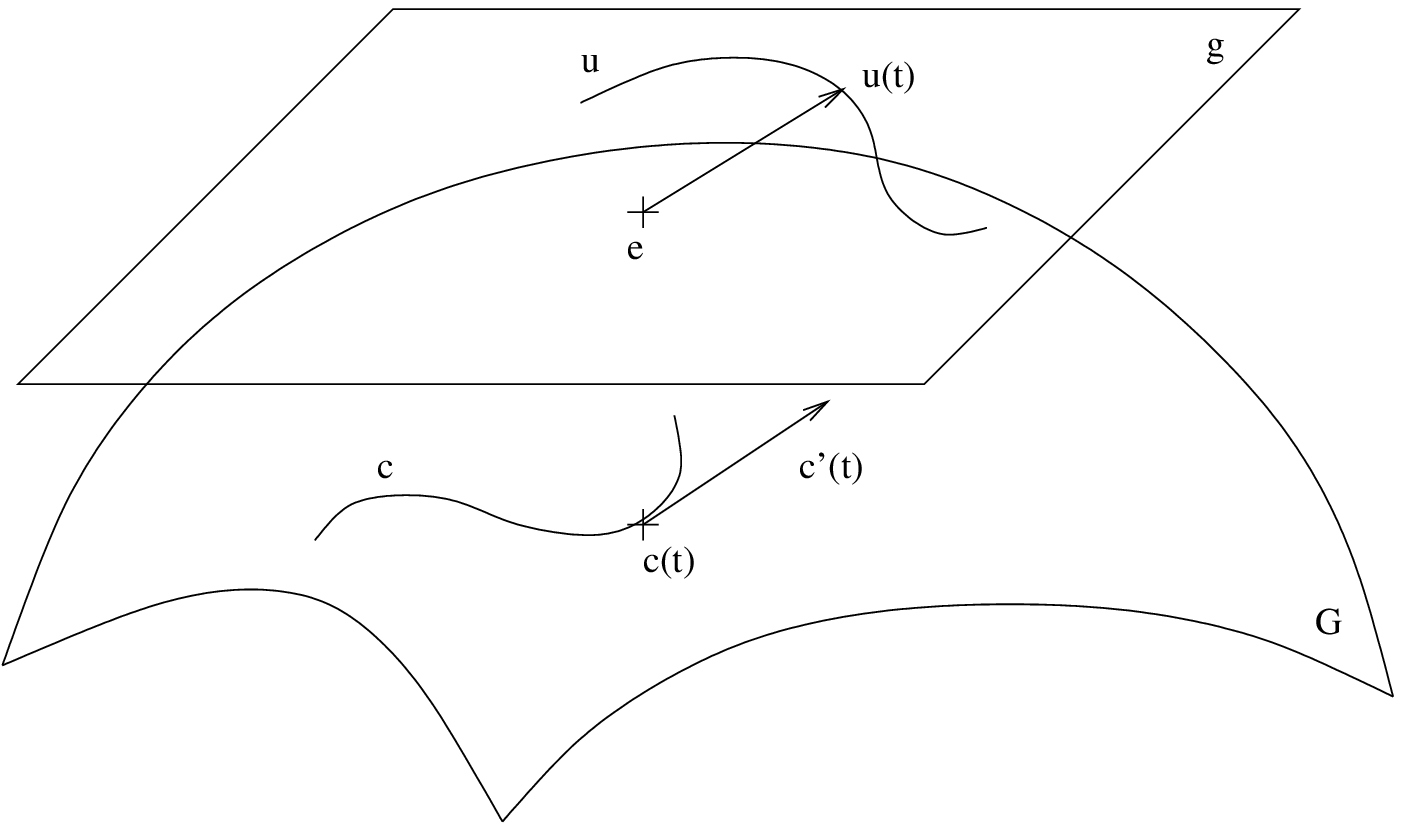}}
\end{figure}

\begin{proof}
We denote the given curve by $c_0$ and its logarithmic derivative by $u_0$.
For any variation with f\/ixed endpoints $c(t,s)\in G$, $t\in [a,b]$, $s\in(-\ep,\ep)$
of the given curve $c_0$, we def\/ine
$u=(\partial_tc)c^{-1}$ and
$v=(\partial_sc)c^{-1}$. In particular $u(\cdot,0)=u_0$, and we denote $v(\cdot,0)$ by $v_0$.

Following \cite{Milnor} we show f\/irst that
\begin{gather}\label{mc}
\partial_t v-\partial_s u=[u,v].
\end{gather}
For each $h\in G$ we consider the map
$F_h(t,s)=(t,s,c(t,s)h)$ for $t\in[a,b]$ and $s\in (-\ep,\ep)$.
The bracket of the following two vector f\/ields on $[a,b]\times (-\ep,\ep)\x G$ vanishes:
\begin{gather*}
(t,s,g)\mapsto\partial_t+u(t,s)g,\qquad (t,s,g)\mapsto \partial_s+v(t,s)g.
\end{gather*}
The reason is  they correspond under the mappings $F_h$, $h\in G$, to the vector f\/ields
$\partial_t$ and $\partial_s$ on $ [a,b]\times (-\ep,\ep)$ (with vanishing bracket).
Hence
$0=[\partial_t+ug,\partial_s+vg]=(\partial_t v)g-(\partial_s u)g-[u,v]g$,
because the bracket of right invariant vector f\/ields corresponds to the opposite bracket
on the Lie algebra $\g$,
so the claim (\ref{mc}) follows.

As in \cite{MR} we compute the derivative of
$E(c)=\frac12\int_a^b\langle u,u\rangle dt$ with respect to $s$,
using the fact that $v(a,s)=v(b,s)=0$.
\begin{gather*}
\partial_s
E(c)
=\int_a^b\langle\partial_su,u\rangle dt\stackrel{(\ref{mc})}{=}
\int_a^b\langle\partial_tv-[u,v],u\rangle dt
=-\int_a^b\langle v,\partial_t u+\ad(u)^\top u\rangle dt.
\end{gather*}
The curve $c_0$ in $G$ is a geodesic if and only if this derivative vanishes at $s=0$
for all variations~$c$ of~$c_0$, hence for all $v_0:[a,b]\to\g$.
This is equivalent to $\frac{d}{dt} u_0=-\ad(u_0)^\top u_0$.
\end{proof}

The Euler equation for a left invariant metric on a Lie group is $\frac{d}{dt} u=\ad(u)^\top u$.
In the case $G=SO(3)$ one obtains the equations of the rigid body.

Denoting by $(\ , \ )$ the pairing between $\g^*$ and $\g$,
the {\it inertia operator} \cite{AK} is def\/ined by
\begin{gather*}
A:\g\to\g^*,\qquad A(X)=\langle X,\cdot\rangle,\qquad \text{\ie} \qquad (A(X),Y)=\langle X,Y\rangle,
\qquad\forall\, X,Y\in\g.
\end{gather*}
It is injective, but not necessarily surjective for inf\/inite dimensional $\g$.
The image of $A$ is called {\it the regular part of the dual}
and is denoted by $\g^*_{\rm reg}$.

Let $\ad^*$ be the coadjoint action of $\g$ on $\g^*$ given by
$(\ad^*(X)m,Y)=(m,-\ad(X)Y)$, for $m\in\g^*$.
The inertia operator relates $\ad(X)^\top$ to the opposite of the coadjoint action of $X$,
\ie
\begin{gather}\label{iden}
\ad^*(X)A(Y)=-A(\ad(X)^\top Y).
\end{gather}
Hence the inertia operator transforms the Euler equation (\ref{euler}) into an equation for $m=A(u)$:
\begin{gather}\label{hami}
\frac{d}{dt}m=\ad^*(u)m,
\end{gather}
result known also as the second Euler theorem.

First Euler theorem states that the solution of (\ref{hami}) with $m(a)=m_0$ is
\[
m(t)=\Ad^*(c(t))m_0,
\]
where $u=\delta^r c$ and $c(a)=e$. Indeed, $\frac{d}{dt}m=\ad^*(\delta^r c)\Ad^*(c) m_0
=\ad^*(u)m$.

\begin{remark}
Equation (\ref{hami}) is a Hamiltonian equation on $\g^*$ with the canonical Poisson bracket
\begin{gather*}
\{f,g\}(m)=\Big(m,\Big[\frac{\de f}{\de m},\frac{\de g}{\de m}\Big]\Big),\qquad f,g\in C^\oo(\g^*)
\end{gather*}
and the Hamiltonian function $h\in C^\oo(\g^*)$,
$h(m)
=\frac12(m,A^{-1}m)=\frac12(m,u)$.
\end{remark}

\begin{remark}
The Euler--Lagrange equation for a right
invariant Lagrangian $L:TG\to\RR$ with value $l:\g\to\RR$ at the identity is:
\begin{gather*}
\frac{d}{dt}\frac{\de l}{\de u}=\ad^*(u)\frac{\de l}{\de u},
\end{gather*}
also called the right Euler--Poincar\'e equation \cite{Poincare,MR0}.
The Hamiltonian form (\ref{hami}) of Euler's equation is obtained for
$l(u)=\frac12\langle u,u\rangle$ since the functional derivative
$\frac{\de l}{\de u}$ is $A(u)$ in this case.
\end{remark}


\section
{Ideal hydrodynamics}\label{3}

Let $G=\Diff_{\mu}(M)$ be the regular Fr\' echet Lie group
of volume preserving dif\/feomorphisms
of a~compact Riemannian manifold $(M,g)$ with induced volume form $\mu$.
Its Lie algebra is $\g=\mathfrak{X}_{\mu}(M)$,
the Lie algebra of divergence free vector f\/ields,
with Lie bracket the opposite of the usual bracket of vector f\/ields
$\ad(X)Y=-[X,Y]$.
We consider the right invariant metric on~$G$ given by the~$L^2$ scalar
product on vector f\/ields
\begin{gather}\label{el2}
\langle X,Y\rangle =\int_M g(X,Y)\mu.
\end{gather}

In the $L^2$ orthogonal decomposition $\X(M)=\X_\mu(M)\oplus\grad(C^\oo(M))$,
we denote by $P$ the projection on $\X_\mu(M)$.
The adjoint of $\ad(X)$
is $\ad(X)^\top Y=P(\nabla_X Y+(\nabla X)^\top Y)$ where
$\nabla$ denotes the Levi-Civita covariant derivative. Indeed,
\begin{gather*}
\langle \ad(X)^\top Y,Z\rangle
=\int_Mg(Y,[Z,X])\mu
=\int_Mg(Y,\nabla_ZX-\nabla_XZ)\mu\\
\phantom{\langle \ad(X)^\top Y,Z\rangle}{}=\int_Mg((\nabla X)^\top Y,Z)\mu+\int_Mg(\nabla_X Y,Z)\mu
=\langle P(\nabla_X Y+(\nabla X)^\top Y),Z\rangle,
\end{gather*}
with $(\nabla X)^\top$ denoting the adjoint of the (1,1)-tensor $\nabla X$
relative to the metric: $g(\nabla_ZX,Y)=g(Z,(\nabla X)^\top Y)$.
In particular $\ad(X)^\top X=P(\nabla_XX)=\nabla_XX+\grad p$,
with $p$ the smooth function uniquely def\/ined up to a constant
by $\De p=\div(\nabla_XX)$.
Now Theorem \ref{theo} assures that
the geodesic equation in $\Diff_\mu(M)$, in terms of the right logarithmic derivative
$u$ of the geodesic, is {\bf Euler's equation for ideal f\/low} with velocity $u$ and pressure $p$
\cite{Moreau,Arnold,EM}:
\begin{gather}\label{ihd}
\partial_tu=-\nabla_uu-\grad p,\qquad\div u=0.
\end{gather}

The geodesic equation (\ref{ihd})
written for the vorticity 2-form $\om=du^\flat$,
$\flat$ denoting the inverse of the Riemannian lift $\sharp$ and $L$ the Lie derivative, is
\begin{gather}\label{vort}
\partial_t\om=-L_u\om,
\end{gather}
because $(\nabla_uu)^\flat=L_uu^\flat-\frac12d(g(u,u))$ and $(\grad p)^\flat=dp$.


\section{Burgers equation}

Let $G=\Diff(S^1)$ be the group of orientation preserving
dif\/feomorphisms of the circle and $\g=\mathfrak{X}(S^1)$
the Lie algebra of vector f\/ields. The Lie bracket is $[X,Y]=X'Y-XY'$,
the negative of the usual bracket on vector f\/ields (vector
f\/ields on the circle are identif\/ied here with their coef\/f\/icient functions in $C^\oo(S^1)$).
We consider the right invariant metric on $G$ given by the $L^2$ scalar
product
$\langle X,Y\rangle =\int_{S^1}XYdx$ on $\g$. The adjoint of $\ad(X)$
is $\ad(X)^\top Y=2X'Y+XY'$, because:
\begin{gather*}
\langle \ad(X)^\top Y,Z\rangle
=\int_{S^1}Y(X'Z-XZ')dx
=\int_{S^1}(X'Y+(XY)')Zdx=\langle 2X'Y+XY',Z\rangle .
\end{gather*}
It follows from Theorem \ref{theo} that
the geodesic equation on $\Diff(S^1)$ in terms of the right logarithmic derivative
$u:I\to C^\oo(S^1)$ is {\bf Burgers inviscid equation} \cite{Burgers}:
\begin{gather}\label{burgers}
\partial_tu=-3uu'.
\end{gather}

The higher dimensional Burgers equation is the {\bf template
matching equation}, used for comparing images via a deformation induced distance.
It is the geodesic equation on $\Diff(M)$,
the dif\/feomorphism group of a compact Riemannian manifold $(M,g)$, for
the right invariant $L^2$ metric \cite{HMA,V3}:
\begin{gather}\label{matching}
\partial_tu=-\nabla_uu-(\div u)u-\tfrac12\grad g(u,u).
\end{gather}
Indeed,
\begin{gather}\label{brasov}
\ad(X)^\top=(\div X)1+\nabla_X+(\nabla X)^\top,\qquad\forall \, X\in\X(M),
\end{gather}
because as in Section \ref{3} we compute
$\langle \ad(X)^\top Y,Z\rangle
=\int_Mg((\nabla X)^\top Y,Z)\mu+\int_Mg(\nabla_X Y,Z)\mu
-\int_ML_Xg(Y,Z)\mu
=\langle (\nabla X)^\top Y+\nabla_XY+(\div X)Y,Z\rangle$
for all vector f\/ields $X$, $Y$, $Z$ on $M$.

In particular for $M=S^1$ and $u$ a curve in $\X(S^1)$, identif\/ied with $C^\oo(S^1)$, $\div u=u'$ and $g(u,u)=u^2$, so
each of the three terms in the right hand side of (\ref{matching})
is $-uu'$ and we recover Burgers equation (\ref{burgers}).


\section{Abelian extensions}

A bilinear skew-symmetric map $\om:\g\x\g\to V$ is a 2-cocycle on the Lie algebra $\g$
with values in the $\g$-module $V$ if it satisf\/ies the condition
\begin{gather*}
\sum_{\rm cycl}\om([X_1,X_2],X_3)=\sum_{\rm cycl}b(X_1)\om(X_2,X_3), \qquad X_1,X_2,X_3\in\g,
\end{gather*}
where $b:\g\to L(V)$ denotes the Lie algebra action on $V$.
It determines an Abelian
Lie algebra extension $\hat\g:=V\rtimes_\om\g$ of $\g$ by the $\g$-module $V$ with Lie bracket
\begin{gather}\label{bra}
[(v_1,X_1),(v_2,X_2)]=(b(X_1)v_2-b(X_2)v_1+\om(X_1,X_2),[X_1,X_2]).
\end{gather}
There is a 1-1 correspondence between
the second Lie algebra cohomology group $H^2(\g,V)$ and equivalence classes of Abelian
Lie algebra extensions $0\to V\to\hat\g\to\g\to 0$.

When $G$ is inf\/inite dimensional, the two obstructions
for the integrability of such an Abelian Lie algebra
extension to a Lie group extension of the connected Lie group $G$ involve $\pi_1(G)$ and~$\pi_2(G)$~\cite{Neeb}. The Lie algebra 2-cocycle $\om$ is integrable if
\begin{itemize}\itemsep=0pt
\item the period group $\Pi_\om\subset V$ (the group of spherical periods
of the equivariant $V$-valued 2-form on $G$ def\/ined by $\om$) is discrete and
\item the f\/lux homomorphism $F_\om:\pi_1(G)\to H^1(\g,V)$ vanishes.
\end{itemize}
Then for any discrete subgroup $\Ga$ of the subspace of $\g$-invariant elements
of $V$ with $\Ga\supseteq\Pi_\om$, there is an Abelian Lie group extension
$1\to T\to\hat G\to G\to 1$ of $G$ by $T=V/\Ga$.

There are two special cases:
\begin{enumerate}\itemsep=0pt
\item {\bf Semidirect product}: $\hat\g=V\rtimes\g$, obtained when $\om=0$.\\
An example is the semidirect product $\g^*\rtimes G$
for the coadjoint $G$-action on $\g^*$, called the magnetic extension in \cite{AK}.
It has the Lie algebra $\g^*\rtimes\g$, a semidirect  product for
the coadjoint $\g$-action $b=\ad^*$ on $\g^*$.

\item {\bf Central extension}: $\hat\g=V\x_\om\g$, obtained when $b=0$.\\
An example is the Virasoro algebra $\RR\x_\om\X(S^1)$,
a central extension of the Lie algebra of vector f\/ields on the circle
given by the Virasoro cocycle
$\om(X,Y)=\int_{S^1}(X'Y''-X''Y')dx$.
It has a corresponding Lie group extension of the group $\Diff(S^1)$ of orientation
preserving dif\/feomorphisms of the circle,
def\/ined by the Bott group cocycle:
\begin{gather}\label{bott}
c(\ph,\ps)=\int_{S^1}\log(\ph'\circ\ps)d\log\ps',\qquad \ph,\ps\in\Diff(S^1).
\end{gather}
\end{enumerate}

An example of a general Abelian Lie algebra extension is $C^\oo(M)\rtimes_\om\X(M)$,
the Abelian extension of the Lie algebra of vector f\/ields on the manifold $M$
with the opposite bracket
by the natural module of smooth functions on $M$,
the Lie algebra action being $b(X)f=-L_Xf$.
The cocycle $\om:\X(M)\x\X(M)\to C^\oo(M)$
is given by a closed dif\/ferential 2-form $\et$ on $M$.
If $\et$ is an~integral form, then there is a principal circle bundle
$P$ over $M$ with curvature $\et$. In this case the group of
equivariant automorphisms of $P$ is a Lie group extension
integrating the Lie algebra cocycle $\om$:
\begin{gather}\label{gauge}
1\to C^\oo(M,\TT)\to\Diff(P)^\TT\to\Diff(M)_{[P]}\to 1.
\end{gather}
Here $C^\oo(M,\TT)$ is the gauge group of $P$ and $\Diff(M)_{[P]}$ is the group of dif\/feomorphisms of $M$ preserving the bundle class $[P]$ under pullbacks (group having the same identity component as $\Diff(M)$).


\section{Geodesic equations on Abelian extensions}\label{6}

Following \cite{V1} we write down the geodesic equations on an Abelian
Lie group extension
$\hat G$ of $G$ with respect to the right invariant metric def\/ined with the
scalar product
\begin{gather}\label{sca}
\langle (v_1,X_1),(v_2,X_2)\rangle _{\hat\g}=
\langle v_1,v_2\rangle _V+\langle X_1,X_2\rangle _\g
\end{gather}
on its Lie algebra $\hat\g=V\rtimes_\om\g$.
Here $\langle \ , \ \rangle _\g$ and $\langle \ , \ \rangle _V$ are scalar products on
$\g$ and $V$. We have to assume the existence of the following maps:
the adjoint $\ad(X)^\top:\g\to\g$ and the adjoint
$b(X)^\top:V\to V$ for any $X\in\g$,
the linear map $h:V\to L_{\rm skew}(\g)$ taking values in the space of
skew-adjoint operators on $\g$, def\/ined by
\begin{gather*}
\langle h(v)X_1,X_2\rangle _\g=\langle \om(X_1,X_2),v\rangle _V,
\end{gather*}
and the bilinear map $l:V\x V\to\g$,
def\/ined by
\begin{gather*}
\langle l(v_1,v_2),X\rangle _\g=\langle b(X)v_1,v_2\rangle _V.
\end{gather*}
The diamond operation $\diamond:V\x V^*\to\g$ in \cite{HMR} corresponds to our map $l$
via $\langle \ , \ \rangle_V$.

\begin{proposition}\label{Abelian}
The geodesic equation on the Abelian extension $\hat G$ for the right invariant metric
defined by the scalar product \eqref{sca} on $\hat\g$, written
for the right logarithmic derivative $(f,u)$, \ie for curves $u$ in $\g$ and $f$ in $V$, is
\begin{gather*}
\frac{d}{dt}u =-\ad(u)^\top u-h(f)u+l(f,f),\\
\frac{d}{dt}f =-b(u)^\top f.
\end{gather*}
\end{proposition}

\begin{proof}
We compute the adjoint of $\ad(v,X)$ in $V\rtimes_\om\g$ \wrt the scalar product (\ref{sca})
\begin{gather*}
\langle\ad(v_1,X_1)^\top (v_2,X_2),(v_3,X_3)\rangle_{\hat\g}
=\langle(v_2,X_2),(b(X_1)v_3-b(X_3)v_1+\om(X_1,X_3),[X_1,X_3]\rangle_{\hat\g}\\
\qquad{} =\langle v_2,b(X_1)v_3\rangle_V+\langle X_2,[X_1,X_3]\rangle_\g
+\langle v_2,\om(X_1,X_3)\rangle_V-\langle v_2,b(X_3)v_1\rangle_V\\
\qquad{} =\langle(b(X_1)^\top v_2,\ad(X_1)^\top X_2+h(v_2)X_1-l(v_1,v_2)),(v_3,X_3)\rangle_{\hat\g}.
\end{gather*}
The result follows now from Euler's equation (\ref{euler}).
\end{proof}

\begin{remark}\label{rema1}
When the scalar product on $V$ is $\g$-invariant,
\ie $\langle b(X)v_1,v_2\rangle_V+\langle v_1,b(X)v_2\rangle_V=0$,
then $l$ is skew-symmetric and the geodesic equation becomes
\begin{gather*}
\frac{d}{dt}u =-\ad(u)^\top u-h(f)u,\\
\frac{d}{dt}f =b(u) f.
\end{gather*}
\end{remark}


\section{Geodesic equations on semidirect products}\label{7}


A special case of Proposition \ref{Abelian}, obtained for $\om=0$, is:

\begin{corollary}\label{semidirect}
The geodesic equation on the semidirect product Lie group
$V\rtimes G$ for the right invariant metric
defined by the scalar product \eqref{sca}, written
for the curve $(f,u)$ in $V\rtimes\g$, is
\begin{gather*}
\frac{d}{dt}u =-\ad(u)^\top u+l(f,f),\\
\frac{d}{dt}f =-b(u)^\top f.
\end{gather*}
It reduces to
\begin{gather*}
\frac{d}{dt}u =-\ad(u)^\top u,\\
\frac{d}{dt}f =b(u) f
\end{gather*}
when the scalar product on $V$ is $\g$-invariant.
\end{corollary}

\subsection*{Passive scalar motion}

The geodesic equation on the semidirect product $C^\oo(M)\rtimes\Diff_\mu(M)$
with $L^2$ right invariant metric, written for the right logarithmic
derivative $(f,u):I\to C^\oo(M)\rtimes\X_\mu(M)$ models {\bf passive scalar
motion} \cite{Hattori}:
\begin{gather}
\partial_tu =-\nabla_uu-\grad p,\nonumber\\
\partial_tf =-df(u).\label{star}
\end{gather}
In this case the $L^2$ scalar product on $C^\oo(M)$ is $\X_\mu(M)$-invariant
and we apply Corollary~\ref{semidirect} to get this geodesic equation.

\section{Magnetohydrodynamics}

Let $A:\g\to\g^*$ be the inertia operator def\/ined by a f\/ixed
scalar product $\langle \ ,\ \rangle $ on $\g$.
The scalar product on the regular dual $\g^*_{\rm reg}=A(\g)$
induced via $A$ by this scalar product in $\g$ is again denoted by $\langle \ , \ \rangle $.
Next we consider the subgroup $\g^*_{\rm reg}\rtimes G$ of the magnetic extension
$\g^*\rtimes G$,
with right invariant metric of type (\ref{sca}) \cite{V6}.

\begin{proposition}\label{corollary}
If the adjoint of $\ad(X)$ exists for any $X\in\g$, then
the geodesic equation on the magnetic extension $\g^*_{\rm reg}\rtimes G$
with right invariant metric,
written for the curve $(A(v),u)$ in $\g^*_{\rm reg}\rtimes\g$
is
\begin{gather*}
\frac{d}{dt}u =-\ad(u)^\top u+\ad(v)^\top v,\\
\frac{d}{dt}v =\ad(u)v.
\end{gather*}
\end{proposition}

\begin{proof}
We have to compute the map $l:\g^*_{\rm reg}\x\g^*_{\rm reg}\to\g$
and the adjoint $b(X)^\top:\g^*_{\rm reg}\to\g^*_{\rm reg}$ for $b=\ad^*$.
We use the fact (\ref{iden}) that the coadjoint action on the image of $A$
comes from the opposite of $\ad(\cdot)^\top$.
Then $l(A(Y_1),A(Y_2))=\ad(Y_2)^\top Y_1$
because
\begin{gather*}
\langle l(A(Y_1),A(Y_2)),X\rangle=\langle\ad^*(X)A(Y_1),A(Y_2)\rangle
=-\langle \ad(X)^\top Y_1,Y_2\rangle=\langle\ad(Y_2)^\top Y_1,X\rangle.
\end{gather*}
Also the adjoint of $b(X)=\ad^*(X)$ exists
and $b(X)^\top A(Y)=-A(\ad(X)Y)$.
The result follows now from Corollary \ref{semidirect}.
\end{proof}

For $G=SO(3)$ and left invariant metric on its magnetic extension $\g^*\rtimes G$
one obtains Kirchhof\/f equations for a rigid body moving in a f\/luid.

Let $G=\Diff_\mu(M)$ be the group of volume preserving dif\/feomorphisms on
a compact mani\-fold $M$ and $\g=\X_\mu(M)$.
The regular part $\g^*_{\rm reg}$ of $\g^*$ is naturally isomorphic to the
quotient space $\Om^1(M)/d\Om^0(M)$ of dif\/ferential 1-forms modulo exact 1-forms,
the pairing being $([\al],X)=\int_M\al(X)\mu$, for $\al\in\Om^1(M)$.
More precisely $A(X)$ is the coset $[X^\flat]$ obtained via the Riemannian metric.

Considering the right invariant $L^2$ metric on the magnetic extension $\g^*_{\rm reg}\rtimes G$
determined by the $L^2$ scalar product (\ref{el2}) on vector f\/ields,
the geodesic equations for the time dependent divergence free vector f\/ields $u$ and $B$ are
(by Proposition \ref{corollary})
\begin{gather*}
\partial_tu =-\nabla_u u+\nabla_B B-\grad p,\\
\partial_tB =-L_uB.
\end{gather*}

We specialize to a three dimensional manifold $M$.
The curl of a vector f\/ield $X$ is the vector f\/ield
def\/ined by the relation $i_{\curl X}\mu=dX^\flat$ and the cross product of two vector
f\/ields $X$ and~$Y$ is the vector f\/ield def\/ined by the relation
$(X\x Y)^\flat=i_Yi_X\mu$.
A short computation gives
$(\curl X\x X)^\flat=i_XdX^\flat=L_XX^\flat-dg(X,X)=(\nabla_XX)^\flat-\frac12dg(X,X)$,
hence $\nabla_XX=\curl X\x X+\frac12\grad g(X,X)$.
The geodesic equations above are in this case the equations of
{\bf ideal magnetohydrodynamics}
with velocity $u$, magnetic f\/ield $B$ and pressure $p$
\cite{VD,MRW}:
\begin{gather*}
\partial_tu =-\nabla_u u+\curl B\x B-\grad p,\\
\partial_tB =-L_uB.
\end{gather*}


\subsection*{Magnetic hydrodynamics with asymmetric stress tensor}

Let $M$ be a 3-dimensional compact parallelizable Riemannian manifold
with induced volume form $\mu$
and let $G=\Diff_\mu(M)$ with $\g=\X_\mu(M)$.
Each vector f\/ield $X$ on $M$ can be identif\/ied with a smooth function in
$C^\oo(M,\RR^3)$, and $j(X)\in C^\oo(M,\gl(3,\RR))$ denotes its Jacobian.
Then
$\omega(X,Y)=[\tr(j(X)dj(Y))]\in\Omega^1(M)/d\Omega^0(M)$ is a Lie algebra 2-cocycle on $\g$
with values in the regular dual $\g^*_{\rm reg}$.

Considering the $L^2$ scalar product on the Abelian extension $\g_{\rm reg}^*\rtimes_\om\g$,
we get the following Euler equation \cite{Billig} for time dependent divergence free vector fields $u$ and $B$:
\begin{gather*}
\partial_tu =-\nabla_u u+\curl B\x B+\tr (j(B)\grad j(u))-\grad p,\\
\partial_tB =-L_uB,
\end{gather*}
modeling {\bf magnetic hydrodynamics with asymmetric stress tensor} $T=j(B)\circ j(u)$.

\section{Geodesic equations on central extensions}

When $V=\RR$ is the trivial $\g$-module, then the Lie algebra action $b$ vanishes and
we get a central extension $\RR\times_\om\g$ def\/ined by the cocycle
$\om:\g\times\g\to\RR$.
A consequence of Proposition \ref{Abelian} is:

\begin{corollary}\label{central}
The geodesic equation on a $1$-dimensional central Lie group extension
$\hat G$ of $G$ with right invariant metric determined by the scalar product
$\langle (a,X),(b,Y)\rangle_{\hat\g}=\langle X,Y\rangle_\g+ab$
on its Lie algebra $\hat\g=\RR\times_\om\g$ is
\begin{gather*}
\frac{d}{dt}u=-\ad(u)^\top u-ak(u),\qquad a\in\RR,
\end{gather*}
where $u$ is a curve in $\g$ and $k\in L_{\rm skew}(\g)$ is defined by the Lie algebra cocycle $\om$
via
\begin{gather*}
\langle k(X),Y\rangle =\om(X,Y),\qquad\forall \, X,Y\in\g.
\end{gather*}
\end{corollary}

\begin{proof}
The central extension is a particular case of an Abelian extension, so
Proposition \ref{Abelian} can be applied. The linear map $h:\RR\to L_{\rm skew}(\g)$
has the form $h(a)X=ak(X)$,
because $\langle h(a)X_1,X_2\rangle_\g=a\om(X_1,X_2)=\langle ak(X_1),X_2\rangle_\g$.
The $\g$-module $\RR$ being trivial, $\frac{d}{dt}a=0$, so $a\in\RR$ is constant.
\end{proof}

\subsection*{KdV equation}

The geodesic equation on the Bott--Virasoro group (\ref{bott})
for the right invariant $L^2$ metric is the {\bf Korteweg--de Vries} equation
\cite{OK}.
In this case the Lie algebra is the central extension of $\g=\X(S^1)$
(identif\/ied with $C^\oo(S^1)$) given by
the Virasoro cocycle $\om(X,Y)=\int_{S^1}(X'Y''-X''Y')dx$.
The computation $\om(X,Y)=-2\int_{S^1}X''Y'dx=2\int_{S^1}X'''Ydx=\langle X''',Y\rangle$
implies $k(X)=2X'''$
and by Corollary \ref{central} the geodesic equation for $u:I\to C^\oo(S^1)$ is the KdV equation:
\begin{gather*}
\partial_tu=-3uu'-2au''',\qquad a\in\RR.
\end{gather*}

\section{Superconductivity equation}\label{sec10}

Given a compact manifold $M$ with volume form $\mu$,
each closed 2-form $\et$ on $M$ def\/ines a Lichnerowicz 2-cocycle $\om_\et$
on the Lie algebra of divergence free vector f\/ields,
\begin{gather*}
\om_\et(X,Y)=\int_M\et(X,Y)\mu.
\end{gather*}

The kernel of the f\/lux homomorphism
\begin{gather*}
\flux_\mu:X\in\X_\mu(M)\mapsto [i_X\mu]\in H^{n-1}(M,\RR)
\end{gather*}
is the Lie algebra $\X_\mu^{\rm ex}(M)$ of exact divergence free vector f\/ields.
On a 2-dimensional manifold it consists of vector f\/ields $X$ possessing
stream functions $f\in C^\oo(M)$,
\ie $i_X\mu=df$ ($X$ is the Hamiltonian vector f\/ield with Hamiltonian function $f$).
On a 3-dimensional manifold
it consists of vector f\/ields $X$ possessing vector potentials $A\in\X(M)$,
\ie $i_X\mu=dA^\flat$ ($X$ is the curl of $A$).

The Lie algebra homomorphism $\flux_\mu$ integrates to the f\/lux homomorphism
(due to Thurston) $\Flux_\mu$ on the identity component of the group
of volume preserving dif\/feomorphisms:
\begin{gather*}
\Flux_\mu:\Diff_\mu(M)_0\to H^{n-1}(M,\RR)/\Ga,\qquad
\Flux_\mu(\ph)=\int_0^1[i_{\de^r\ph(t)}\mu]dt\mod\Ga,
\end{gather*}
where $\ph(t)$ is any volume preserving dif\/feotopy from the identity on $M$ to $\ph$
and $\Ga$ a discrete subgroup of $H^{n-1}(M,\RR)$.
The kernel of $\Flux_\mu$ is, by def\/inition, the Lie group $\Diff_\mu^{\rm ex}(M)$
of {\it exact volume preserving diffeomorphisms}.
It coincides with $\Diff_\mu(M)_0$ if and only if $H^{n-1}(M,\RR)=0$.

For $\et$ integral, the Lichnerowicz cocycle is integrable to $\Diff_\mu^{\rm ex}(M)$
\cite{Ismagilov}.
When $M$ is 3-dimensional, there exists a vector f\/ield $B$ on $M$ def\/ined
with $\et=-i_B\mu$. The 2-form $\et$ is
closed if and only if $B$ is divergence free.
The integrality condition of $\et$ expresses as
$\int_S(B\cdot n)d\si\in\ZZ$ on every closed surface $S\subset M$.

The {\bf superconductivity equation} models the
motion of a high density electronic gas in a~magnetic f\/ield $B$
with velocity $u$:
\begin{gather}\label{superconductivity}
\partial_tu=-\nabla_uu-au\x B-\grad p,\qquad a\in\RR.
\end{gather}
It is the geodesic equation on
a central extension of the group of volume preserving dif\/feomorphisms
for the right invariant $L^2$ metric \cite{Zeitlin2,V1},
when $M$ is simply connected.

Indeed,
\begin{gather*}
\om_\et(X,Y)=\int_M\et(X,Y)\mu=-\int_M\mu(B,X,Y)\mu=\int_M g(X\x B,Y)\mu
=\langle P(X\x B),Y\rangle
\end{gather*}
hence the map $k\in L_{\rm skew}(\g)$ determined by the Lichnerowicz cocycle $\om_\et$ is
$k(X)=P(X\x B)$, with $P$ denoting the orthogonal projection
on the space of divergence free vector f\/ields. Now we apply Corollary \ref{central}.



\section[Charged ideal fluid]{Charged ideal f\/luid}\label{11}

Let $M$ be an $n$-dimensional Riemannian manifold with
Levi-Civita connection $\nabla$ and volume form $\mu$, and
$\et$ a closed integral dif\/ferential two-form.
Let $B$ be an $(n-2)$ vector f\/ield on $M$
(i.e. $B\in C^\oo(\wedge^{n-2}TM)$)
such that $\et=(-1)^{n-2}i_B\mu$ is a closed
two-form. The cross product of a~vector f\/ield $X$ with $B$ is
the vector f\/ield
$X\x B=(i_{X\wedge B}\mu)^\sharp=(i_X\et)^\sharp$,
$\sharp$ denoting the Riemannian lift.
When $M$ is 3-dimensional, then $B$ is a divergence free vector f\/ield
with $\et=-i_B\mu$
and $\times$ is the cross product of vector f\/ields.

 From the integrality of $\et$ follows the existence of a~principal $\TT$-bundle $\pi:P\to M$
with a~principal connection 1-form $\al$ on $P$ having curvature $\et$.
The associated Kaluza--Klein metric~$\kappa$ on~$P$, def\/ined at a point~$x\in P$
by
\begin{gather*}
\kappa_x(\tilde X,\tilde Y)=g_{\pi(x)}
(T_x\pi.\tilde X,T_x\pi.\tilde Y)+\al_x(\tilde X)\al_x(\tilde Y),\qquad
\tilde X,\tilde Y\in T_xP
\end{gather*}
determines the volume form $\tilde\mu=\pi^*\mu\wedge\al$ on $P$.

The group $\Diff_{\tilde\mu}(P)^\TT$ of volume preserving automorphisms of the
principal bundle $P$
is an Abelian Lie group extension of $\Diff_\mu(M)_{[P]}$,
the group of volume preserving dif\/feomorphisms
preserving the bundle class $[P]$,
by the gauge group $C^\oo(M,\TT)$
(an extension contained in (\ref{gauge})).
The corresponding Abelian Lie algebra extension
\begin{gather*}
0\to C^\oo(M)\to\X_{\tilde\mu}(P)^\TT\to\X_\mu(M)\to 0
\end{gather*}
is described again by the Lie algebra cocycle $\om:\X_\mu(M)\x\X_\mu(M)\to C^\oo(M)$
given by $\et$.

The Kaluza--Klein metric on $P$ determines a right invariant $L^2$ metric on the group of volu\-me preserving automorphisms of the principal $\TT$-bundle $P$.
The geodesic equation written in terms
of the right logarithmic derivative $(\rho,u)$,
with $\rh$ a time dependent function and $u$ a time dependent
divergence free vector f\/ield on $M$, is:
\begin{gather*}
\partial_tu =-\nabla_uu-\rh u\x B-\grad p,\nonumber\\
\partial_t\rh =-d\rh(u).
\end{gather*}
It models the motion of a {\bf charged ideal f\/luid}
with velocity $u$, pressure $p$ and charge density~$\rh$
in a f\/ixed magnetic f\/ield $B$ \cite{V1}.

Indeed, the connection $\al$ def\/ines a horizontal lift
and identifying the pair $(f,X)$, $f\in C^\oo(M)$, $X\in\X_\mu(M)$
with the sum of the horizontal lift of $X$
and the vertical vector f\/ield given by $f$,
we get an isomorphism between the Abelian Lie algebra
extension $C^\oo(M)\rtimes_\om\X_\mu(M)$
and the Lie algebra $\X_{\tilde\mu}(P)^\TT$ of invariant divergence free vector f\/ields
on $P$. Under this isomorphism the $L^2$ metric
def\/ined by the Kaluza--Klein metric $\kappa$ is
$\langle(f_1,X_1),(f_2,X_2)\rangle=\int_M(g(X_1,X_2)+f_1f_2)\mu$.
The $L^2$ scalar product on functions is $\X_\mu(M)$ invariant, \ie $b(X)$ is skew-adjoint.
The mapping $h:C^\oo(M)\to L_{\rm skew}(\X_\mu(M))$ is $h(f)X=P(fX\x B)$ because:
\begin{gather*}
\langle h(f)X,Y\rangle=\langle\et(X,Y),f\rangle=\int_Mf(i_{X}\et)(Y)\mu
=\int_Mfg(X\times B,Y)\mu=\langle P(fX\times B),Y\rangle,
\end{gather*}
where $P$ denotes the
orthogonal projection on the space of divergence free
vector f\/ields on $M$.
The result follows from Remark \ref{rema1},
knowing that $\ad(X)^\top X=P(\nabla_XX)$.


\section{Geodesics on general extensions}

A general extension of Lie algebras is an exact sequence of Lie algebras
\begin{gather}\label{general}
0\to\h\to\hat\g\to\g\to 0.
\end{gather}
A section $s:\g\to\hat\g$ (\ie a right inverse to the projection $\hat\g\to\g$)
induces the following mappings~\cite{AMR}:
\begin{gather*}
 b: \ \ \g\to\Der(\h),\qquad b(X)f=[s(X),f],\\
 \om:\ \ \g\x\g\to\h,\qquad \om(X_1,X_2)=[s(X_1),s(X_2)]-s([X_1,X_2])
\end{gather*}
with properties:
\begin{gather*}
[b(X_1),b(X_2)]-b([X_1,X_2])=\ad(\om(X_1,X_2)),\\
\sum_{\rm cycl}\om([X_1,X_2],X_3)=\sum_{\rm cycl}b(X_1)\om(X_2,X_3).
\end{gather*}

The Lie algebra structure on the extension $\hat\g$, identif\/ied as a vector space with
$\h\oplus\g$ via the section $s$, can be expressed in terms of $b$ and $\om$:
\begin{gather*}
[(f_1,X_1),(f_2,X_2)]=([f_1,f_2]+b(X_1)f_2-b(X_2)f_1+\om(X_1,X_2),[X_1,X_2]).
\end{gather*}
In particular for $\h$ an Abelian Lie algebra this is the Lie bracket (\ref{bra})
on an Abelian Lie algebra extension.

We consider scalar products $\langle\ ,\ \rangle_\g$ on $\g$ and
$\langle \ , \ \rangle_\h$ on $\h$ and, as in Section~\ref{6},
we impose the existence of several
maps: $\ad(X)^\top:\g\to\g$ for any $X\in\g$,
$\ad(f)^\top:\h\to\h$ for any $f\in\h$,
$b(X)^\top:\h\to\h$ for any $X\in\g$, as well as the linear map
$h:\h\to L_{\rm skew}(\g)$ def\/ined by
\begin{gather*}
\langle h(f)X_1,X_2\rangle _\g=\langle \om(X_1,X_2),f\rangle _\h,
\end{gather*}
and the bilinear map $l:\h\x\h\to\g$,
def\/ined by
\begin{gather*}
\langle l(f_1,f_2),X\rangle _\g=\langle b(X)f_1,f_2\rangle _\h.
\end{gather*}

A result similar to Proposition \ref{Abelian} is:

\begin{proposition}\label{fgb}
The geodesic equation on the Lie group extension $\hat G$
of $G$ by $H$ integra\-ting~\eqref{general}, with right invariant metric
determined by the scalar product
\begin{gather*}
\langle (f_1,X_1),(f_2,X_2)\rangle _{\hat\g}=
\langle f_1,f_2\rangle _\h+\langle X_1,X_2\rangle _\g,
\end{gather*}
written in terms of the right logarithmic derivative $(\rh,u)$ is:
\begin{gather*}
\frac{d}{dt}u =-\ad(u)^\top u-h(\rho)u+l(\rh,\rh),\\
\frac{d}{dt}\rho =-\ad(\rh)^\top\rho-b(u)^\top \rho.
\end{gather*}
\end{proposition}


\section[Ideal fluid in a fixed Yang-Mills field]{Ideal f\/luid in a f\/ixed Yang--Mills f\/ield}

Let $\pi:P\to M$ be a principal $G$-bundle with principal action $\si:P\x G\to P$
and let $\Ad P=P\x_G\g$ be its adjoint bundle.
The space $\Om^k(M,\Ad P)$ of dif\/ferential forms with
values in $\Ad P$ is identif\/ied with the space $\Om^k_{\rm hor}(P,\g)^G$
of $G$-equivariant horizontal forms on $P$. In particular $C^\oo(M,\Ad P)=C^\oo(P,\g)^G$.

We consider a principal connection 1-form $\al\in\Om^1(P,\g)^G$ on $P$.
Its curvature $\et=d\al+\frac12[\al,\al]$.
is an equivariant horizontal 2-form $\et\in\Om^2_{\rm hor}(P,\g)^G$,
hence it can be viewed as a 2-form on $M$ with values in $\Ad P$.
The covariant exterior derivative on $\g$-valued dif\/ferential forms on $P$ is
$d^\al=\chi^*\circ d$,
with $\chi:\X(P)\to\X(P)$ denoting the horizontal projection, and it induces a map
$d^\al:\Om^k(M,\Ad P)\to\Om^{k+1}(M,\Ad P)$.

Let $g$ be a Riemannian metric on $M$ and $\ga$ a $G$-invariant scalar product on $\g$.
These data, together with the connection $\al$, def\/ine a Kaluza--Klein metric on $P$:
\begin{gather*}
\kappa_x(\tilde X,\tilde Y)=g_{\pi(x)}
(T_x\pi.\tilde X,T_x\pi.\tilde Y)+\ga(\al_x(\tilde X),\al_x(\tilde Y)),\qquad
\tilde X,\tilde Y\in T_xP.
\end{gather*}
The canonically induced volume form on $P$ is $\tilde\mu=\pi^*\mu\wedge\al^*\det_\ga$,
where $\mu$ is the canonical volume form on $M$ induced by the Riemannian metric $g$
and $\al^*\det_\ga$ is the pullback by $\al:TP\to\g$ of the determinant
$\det_\ga\in\wedge^{\dim\g}\g^*$
induced by the scalar product $\ga$ on $\g$.

The gauge group of the principal bundle is identif\/ied with $C^\oo(P,G)^G$,
the group of $G$-equivariant functions from $P$ to $G$, with $G$ acting on itself
by conjugation. The group of automorphisms of $P$, \ie the group of $G$-equivariant
dif\/feomorphisms of $P$, is an extension of $\Diff(M)_{[P]}$, the group of dif\/feomorphisms of $M$
preserving the bundle class $[P]$, by the gauge group. This is the analogue of (\ref{gauge})
for non-commutative structure group. Restricting to volume preserving
dif\/feomorphisms, we get the exact sequence:
\begin{gather*}
1\to C^\oo(P,G)^G\stackrel{\si}{\to}\Diff_{\tilde\mu}(P)^G\to\Diff_\mu(M)_{[P]}\to 1.
\end{gather*}

On the Lie algebra level the exact sequence is
\begin{gather*}
0\to C^\oo(P,\g)^G\stackrel{\dot\si}{\to}\X_{\tilde\mu}(P)^G\to\X_\mu(M)\to 0.
\end{gather*}
The horizontal lift provides a linear section $:\X_\mu(M)\to\X_{\tilde\mu}(P)^G$,
thus identifying
the pair $(f,X)\in C^\oo(P,\g)^G\oplus\X_\mu(M)$ with
$\tilde X=\dot\si(f)+X^{\rm hor}\in\X_{\tilde\mu}(P)^G$.
With this identif\/ication,
the~$L^2$ metric on $\X_{\tilde\mu}(P)^G$ given by the Kaluza--Klein metric
can be written as
\begin{gather*}
\int_P\ka(f_1, X_1),(f_2,X_2))\tilde\mu=\int_Mg(X_1,X_2)\mu+\int_P\ga(f_1,f_2)\tilde\mu.
\end{gather*}

A particular case of a result in \cite{GR3} is the fact that the geodesic equation
on the group $\Diff_{\tilde\mu}(P)^G$ of volume preserving automorphisms of $P$
with right invariant $L^2$ metric gives the equations of motion of
an {\bf ideal f\/luid moving in a f\/ixed Yang--Mills f\/ield}.
Written for the right logarithmic derivative
$(\rh,u):I\to C^\oo(P,\g)^G\oplus\X_\mu(M)$, these are:
\begin{gather}
\partial_tu =-\nabla_uu-\ga(\rh,i_u\et)^\sharp-\grad p,\nonumber\\
\partial_t\rh =-d^\al\rh(u).\label{alfa}
\end{gather}
Here $u$ denotes the Eulerian velocity,
$\rh$, viewed as a time dependent section of $\Ad P$, denotes the magnetic charge,
$\et$, viewed as a 2-form on $M$ with values in $\Ad P$,
denotes the f\/ixed Yang--Mills f\/ield. The scalar product $\ga$ being $G$-invariant,
can be viewed as a bundle metric on $\Ad P$.

This result follows from Proposition \ref{fgb}. Indeed, in this particular case
the cocycle is $\om=\et$ and the Lie algebra action is $b(X)f=-df.X^{\rm hor}=-d^\al f.X$, hence
$b(X)^\top f=-b(X)f$, $l$ is skew-symmetric and $h(f)X=P(\ga(f,i_X\et)^\sharp)$. Moreover $\ad(X)^\top X=P\nabla_XX$
with $P$ the projection on divergence free
vector f\/ields and $\ad(f)^\top f=[f,f]=0$,
so~(\ref{alfa}) follows.

The equations of a charged ideal f\/luid from Section \ref{11} are obtained for the structure group~$G$ equal to the torus $\TT$.


\section{Totally geodesic subgroups}\label{14}

Let $G$ be a Lie group with right invariant Riemannian metric.
A Lie subgroup $H\subseteq G$ is totally
geodesic if any geodesic $c:[a,b]\to G$ with $c(a)=e$ and $c'(a)\in\h$,
the Lie algebra of $H$, stays in~$H$.

From the Euler equation (\ref{euler}) we see that this is the case if $\ad(X)^\top X\in\h$ for all
$X\in\h$. If there is a geodesic in $G$ in any direction
of $\h$, then this condition is necessary and suf\/f\/icient,
so we give the following def\/inition: the Lie subalgebra
$\h$ is called {\it totally geodesic} in $\g$ if
$\ad(X)^\top X\in\h$ for all $X\in\h$.

\begin{remark}\label{cap}
Given two totally geodesic Lie subalgebras $\h$ and $\kkk$ of the Lie algebra $\g$,
the intersection $\h\cap\kkk$ is totally geodesic in $\g$, but also in $\h$ and in $\kkk$.
\end{remark}


\subsection*{Ideal f\/luid}

The ideal f\/luid f\/low (\ref{ihd}) on $M$
preserves the property of having a stream function (if $M$ two dimensional),
\resp a vector potential (if $M$ three dimensional)
if and only if $\Diff_\mu^{\rm ex}(M)$ is a totally geodesic subgroup of $\Diff_\mu(M)$
for the right invariant $L^2$ metric.
This means $P(\nabla_XX)\in\X_\mu^{\rm ex}(M)$ for all $X\in\X_\mu^{\rm ex}(M)$.

\begin{theorem}[\cite{HTV}]\label{stefan}
The only Riemannian manifolds $M$ with the property that
$\Diff_\mu^{\rm ex}(M)$
is a~totally geodesic subgroup of $\Diff_\mu(M)$ with the right invariant $L^2$ metric
are twisted products $M=\RR^k\x_\La F$ of a flat torus $\TT^k=\RR^k/\Lambda$ and a connected oriented Riemannian manifold $F$ with $H^1(F,\RR)=0$.
\end{theorem}

In particular the ideal f\/luid f\/low on the 2-torus preserves the property of
having a stream function~\cite{AK} and the ideal f\/luid f\/low
on the 3-torus preserves the property of having a vector potential.

\subsection*{Superconductivity}

Given a compact Riemannian manifold $M$,
from the Hodge decomposition follows that
$\X_\mu(M)=\X_\mu^{\rm ex}(M)\oplus\X_{\rm harm}(M)$.
On a f\/lat torus the harmonic vector f\/ields are those with all components constant.

In the setting of Section \ref{sec10}, the next proposition determines when is $\RR\rtimes_{\om_\et}\X_\mu^{\rm ex}(M)$
totally geodesic in $\RR\rtimes_{\om_\et}\X_\mu(M)$, for $M=\TT^3$ and $\et=-i_B\mu$.

\begin{proposition}[\cite{V2}]
The superconductivity equation \eqref{superconductivity}
on the $3$-torus preserves the pro\-perty of having a vector potential
if and only if the three components of the magnetic field $B$ are constant.
\end{proposition}

\begin{proof}
Any exact divergence free vector f\/ield $X$ on the 3-torus admits a potential 1-form $\al$
with $i_X\mu=d\al$, hence
$\int_{\TT^3}g(X\x B,Y)\mu=\int_{\TT^3}i_Yi_B\mu\wedge i_X\mu
=\int_{\TT^3}i_{[Y,B]}\mu\wedge\al$.
Then the totally geodesicity condition which, in this case, says that
$P(X\x B)$ is exact divergence free for all $X$ exact divergence free,
is equivalent to $[Y,B]=0$ for all harmonic vector f\/ields $Y$.
This is further equivalent to the fact that
the three components of the magnetic f\/ield $B$ are constant.
\end{proof}

\subsection*{Passive scalar motion}\label{11.3}

On the trivial principal $\TT$ bundle $P=M\x\TT$ we consider the volume form
$\tilde\mu=\mu\wedge d\th$. Noticing that $i_{(f,X)}\tilde\mu=i_X\mu\wedge d\th+f\mu$,
we get the Lie algebra isomorphisms
$\X_{\tilde\mu}(M\x\TT)^\TT\cong C^\oo(M)\rtimes\X_\mu(M)$
and $\X_{\tilde\mu}^{\rm ex}(M\x\TT)^\TT\cong C_0^\oo(M)\rtimes\X_\mu^{\rm ex}(M)$,
where $C_0^\oo(M)$ is the subspace
of functions with vanishing integral.

From \cite{V3} we know that the group of equivariant volume preserving dif\/feomorphisms
is totally geodesic in the group of volume preserving dif\/feomorphisms
and from Theorem \ref{stefan} we know that the group of exact volume preserving
dif\/feomorphisms of a torus
is totally geodesic in the group of volume preserving dif\/feomorphisms,
hence by Remark \ref{cap} we obtain that for $M=\TT^2$ the subgroup
$\Diff_{\tilde\mu}^{\rm ex}(M\x\TT)^\TT$
is totally geodesic in $\Diff_{\tilde\mu}(M\x\TT)^\TT$.
This means that $C_0^\oo(M)\rtimes\X_\mu^{\rm ex}(M)$
is totally geodesic in $C^\oo(M)\rtimes\X_\mu(M)$ for $M=\TT^2$.
In other words equation (\ref{star}), describing passive scalar motion,
preserves the property of
having a stream function if $f$ has zero integral at the initial moment.
Moreover, $f$ will have zero integral at any moment.



\section{Quasigeostrophic motion}

Given a closed 1-form $\al$ on the compact symplectic manifold $(M,\si)$,
the Roger cocycle on the Lie algebra $\X_\si^{\rm ex}(M)$ of Hamiltonian vector f\/ields on $M$
is \cite{Roger}
\begin{gather*}
\om_\al(H_f,H_g)=\int_Mf\al(H_g)\si^n.
\end{gather*}
Here $f$ and $g$ are Hamiltonian functions with zero integral
for the Hamiltonian vector f\/ields $H_f$ and $H_g$.
The integrability of the 2-cocycle $\omega_\alpha$ to a central extension of
the group of Hamiltonian dif\/feomorphisms is an open problem.
Partial results are given in \cite{Ismagilov2}.

For $M=\TT^2$ the cocycle $\om_\al$ can be extended to a cocycle
on the Lie algebra of symplectic vector f\/ields $\X_\si(\TT^2)$ by
$\om_\al(\partial_x,\partial_y)
=\om_\al(\partial_x,H_f)=\om_\al(\partial_y,H_f)=0$ \cite{Kirillov}.
The extendability of $\om_\al$ to $\X_\si(M)$ for $M$ an arbitrary symplectic manifold
is studied in \cite{V4}.
To a divergence free vector f\/ield $X$ on the 2-torus one can assign
a smooth function $\ps_X$ on the 2-torus uniquely determined by $X$ through
$d\ps _X=i_X\si-\langle i_X\si\rangle$ and $\int_{\TT^2}\ps _X\si=0$.
Here $\langle \ \rangle$ denotes the average of a 1-form on the torus:
$\langle adx+bdy\rangle=(\int_{\TT^2}a\si)dx+(\int_{\TT^2}b\si)dy$.
In particular $\ps_{H_f}=f$ whenever $f$ has zero integral.

\begin{proposition}[\cite{V5}]
The Euler equation for the $L^2$ scalar product on $\RR\x_{\om_\al}\X_\si(\TT^2)$ is
\begin{gather}\label{sharp}
\partial_tu=-\nabla_uu-\ps _u\al^\sharp-\grad p,
\end{gather}
where the function $\ps _u$ is uniquely determined by $u$ through
$d\ps _u=i_u\si-\langle i_u\si\rangle$ and $\int_{\TT^2}\ps _u\si=0$.
\end{proposition}

\begin{proof}
To apply Corollary \ref{central} we compute the map $k$
corresponding to the cocycle $\om_\al$.
Using the fact that $\om_\al(\partial_x,X)=\om_\al(\partial_y,X)=0$
for all $X\in\X_\si(\TT^2)$, we get
\begin{gather*}
\om_\al(u,X)=\om_\al(H_{\ps _u},X)=\int_{\TT^2}\ps _u\al(X)\si
=\int_{\TT^2}g(\ps _u\al^\sharp,X)\si=\langle P(\ps _u\al^\sharp),X\rangle,
\end{gather*}
hence $k(u)=P(\ps _u\al^\sharp)$.
Knowing also that $\ad(u)^\top u=P(\nabla_uu)$, we get (\ref{sharp}) as the Euler equation
for $a=1$.
\end{proof}

\begin{proposition}[\cite{V5}]\label{totham}
If the two components of the $1$-form $\al$ on $\TT^2$ are constant,
then equa\-tion~\eqref{sharp} preserves the property of having a stream function,
\ie $\RR\x_{\om_\al}\X_\si^{\rm ex}(\TT^2)$ is totally geodesic in
$\RR\x_{\om_\al}\X_\si(\TT^2)$.
In this case the restriction of \eqref{sharp} to Hamiltonian vector fields is
\begin{gather}\label{hfhf}
\partial_tH_\ps =-\nabla_{H_\ps }H_\ps -\ps \al^\sharp-\grad p.
\end{gather}
\end{proposition}

\begin{proof}
By Theorem \ref{stefan}
on the 2-torus $P(\nabla_XX)$ is Hamiltonian for $X$ Hamiltonian,
hence the totally geodesicity condition in this case
is equivalent to the fact that $P(\ps _X\al^\sharp)$
is Hamiltonian for $X$ Hamiltonian. By Hodge decomposition this means $\ps _X\al^\sharp$
is orthogonal to the space of harmonic vector f\/ields, so
\begin{gather*}
\langle P(\ps _X\al^\sharp),Y\rangle
=\int_{\TT^2}g(\ps _X\al^\sharp,Y)\si=\int_{\TT^2}\al(Y)\ps _X\si=0,\qquad\forall\;  Y\text{ harmonic}.
\end{gather*}
On the torus the harmonic vector f\/ields $Y$ are the vector f\/ields with constant components
and the functions $\ps _X$ have vanishing integral by def\/inition,
so the expression above vanishes for all constant vector f\/ields $Y$
if the 1-form $\al$ has constant coef\/f\/icients.
\end{proof}

On the 2-torus with $\si=dx\wedge dy$ and
$u=H_\ps $, the vorticity 2-form is $du^\flat=d(H_\ps)^\flat=(\De\ps)\si$,
hence $\om=\De\ps$ is the vorticity function.
Since $L_u(du^\flat)=L_{H_\ps }(\om\si)=(L_{H_\ps }\om)\si=\{\om,\ps\}\si$,
the vorticity equation (\ref{vort}) written for the vorticity function $\om$ becomes
\begin{gather*}
\partial_t\om=-\{\om,\ps\}.
\end{gather*}
For $\al=\be dy$, $\be\in\RR$, we have
$d(\ps\al^\sharp)^\flat=d\ps\wedge\al=(\be\partial_x\ps)\si$.
Hence the Euler equation (\ref{hfhf})
written for the vorticity function $\om=\De\ps$ with $\ps$ the stream function
of $u$, is the equation for {\bf quasigeostrophic motion in $\be$-plane approximation}
\cite{ZP,HZ}
\begin{gather*}
\partial_t\om=-\{\om,\ps\}-\be\partial_x\ps,
\end{gather*}
with $\be$ the gradient of the Coriolis parameter.


\section{Central extensions of semidirect products}\label{13}

Let $\g$ be a Lie algebra with scalar product $\langle \ , \ \rangle_\g$
and $V$ a $\g$-module with $\g$-action $b$ and $\g$-invariant
scalar product $\langle \ , \ \rangle_V$.
Each Lie algebra 1-cocycle $\al\in Z^1(\g,V)$ (\ie a linear map $\al:\g\to V$
which satisf\/ies
$\al([X_1,X_2])=b(X_1)\al(X_2)-b(X_2)\al(X_1)$)
def\/ines a 2-cocycle $\om$ on the semidirect product $V\rtimes\g$ \cite{OR}:
\begin{gather}\label{or}
\om((v_1,X_1),(v_2,X_2))=\langle\al(X_1),v_2\rangle_V-\langle\al(X_2),v_1\rangle_V.
\end{gather}

\begin{proposition}\label{p6}
The Euler equation on the central extension $(\g\ltimes V)\x_\om\RR$
with respect to the scalar product $\langle \ , \ \rangle_\g+\langle \ , \ \rangle_V$,
written for curves $u$ in $\g$ and $f$ in $V$, is
\begin{gather*}
\frac{d}{dt} u =-\ad(u)^\top u+a\al^\top(f),\nonumber\\
\frac{d}{dt} f =b(u)f-a\al(u), \qquad a\in\RR,
\end{gather*}
where $\al^\top:V\to\g$ is the adjoint of $\al:\g\to V$.
\end{proposition}

\begin{proof}
The map $k\in L_{\rm skew}(V\rtimes\g)$ def\/ined by $\om$
is $k(v,X)=(\al(X),-\al^\top(v))$ because
\begin{gather*}
\om((v_1,X_1),(v_2,X_2))=\langle\al(X_1),v_2\rangle_V-\langle\al^\top(v_1),X_2\rangle_\g
=\langle(\al(X_1),-\al^\top(v_1)),(v_2,X_2)\rangle_{V\rtimes\g}.
\end{gather*}
The result follows from Corollaries
\ref{semidirect} and \ref{central}.
\end{proof}

\begin{remark}\label{ov}
More generally, a 1-cocycle $\al$ on $\g$ with values in the dual $\g$-module $V^*$
def\/ines a~2-cocycle on $V\rtimes\g$ by
\begin{gather*}\label{osienko}
\om((v_1,X_1),(v_2,X_2))=(\al(X_1),v_2)-(\al(X_2),v_1),
\end{gather*}
where $( \ , \ )$ denotes the pairing between $V^*$ and $V$.
\end{remark}


\section[Stratified fluid]{Stratif\/ied f\/luid}

Let $M$ be a compact Riemannian manifold with induced volume form $\mu$.
Let $\al$ be a closed 1-form on $M$. Then $\al:\X(M)\to C^\oo(M)$
is a Lie algebra 1-cocycle with values in the canonical $\X(M)$-module $C^\oo(M)$.
The $L^2$ scalar product is $\X_\mu(M)$-invariant,
so $\al$ def\/ines by (\ref{or}) a~2-cocycle~$\om$ on
the semidirect product Lie algebra $C^\oo(M)\rtimes\X_\mu(M)$:
\begin{gather}\label{oror}
\om((f_1,X_1),(f_2,X_2))=\int_Mf_2\al(X_1)\mu-\int_Mf_1\al(X_2)\mu.
\end{gather}

\begin{proposition}
The Euler equation on $(C^\oo(M)\rtimes\X_\mu(M))\x_\om\RR$ with $L^2$ scalar product is
\begin{gather}
\partial_tu =-\nabla_uu+af\al^\sharp-\grad p,\nonumber\\
\partial_tf =-L_uf-a\al(u),\qquad a\in\RR,\label{our}
\end{gather}
with $\nabla$ the Levi-Civita covariant derivative and $\sharp$ the Riemannian
lift.
\end{proposition}

\begin{proof}
We apply Proposition~\ref{p6} for $\g=\X_\mu(M)$ and $V=C^\oo(M)$.
In this case $b(X)f=-L_Xf$ and $\ad(X)^\top X=P\nabla_XX$. We compute
$\langle\al^\top(f),X\rangle_\g=\int_Mf\al(X)\mu=\int_Mg(f\al^\sharp,X)\mu=
\langle P(f\al^\sharp),X\rangle_\g$, for all $X\in\X_\mu(M)$, hence $\al^\top(f)=
P(f\al^\sharp)$.
\end{proof}

Because $d(f\al^\sharp)^\flat=df\wedge\al$,
the equation (\ref{our}) written for vorticity 2-form $\om=du^\flat$ becomes
\begin{gather*}
\partial_t \om =-L_u\om+df\wedge\al,\\
\partial_t f =-L_uf-a\al(u).
\end{gather*}

\begin{proposition}[\cite{V2}]\label{totgeod}
Given a $2$-cocycle $\om$ determined via \eqref{oror} by the
constant $1$-form $\al$ on the torus $M\!=\!\TT^2$, we have
$(\X_\mu^{\rm ex}(M)\ltimes C^\oo_0(M))\x_\om\RR$ is totally
geodesic in \mbox{$(\X_\mu(M){\ltimes} C^\oo(M)){\x_\om}\RR$},
where $C_0^\oo(M)$ is the subspace
of functions with vanishing integral.
\end{proposition}

\begin{proof}
We know from Section~\ref{14} that for $M=\TT^2$,
$\X_\mu^{\rm ex}(M)\ltimes C^\oo_0(M)$ is totally
geodesic in $\X_\mu(M)\ltimes C^\oo(M)$.
But $\al(u)$ has zero integral for $u$ exact divergence free,
$\al$ being closed. We have to make sure that
$f\al^\sharp$ is orthogonal to the space of harmonic vector f\/ields
for all $f$ with zero integral, \ie for all functions $f$
such that $f\mu$ is exact ($f\mu=d\nu$).
But $\int_Mg(f\al^\sharp,Y)\mu
=\int_M\al(Y)d\nu=-\int_ML_Y\al\wedge\nu=0$
because $L_Y\al=0$ for all harmonic vector
f\/ields $Y$ on the 2-torus, $\al$ being a constant 1-form.
\end{proof}

Hence on the 2-torus,
for constant $\al$ and initial conditions
$u_0$ Hamiltonian vector f\/ield and~$f_0$ function with zero integral,
$u$ will be Hamiltonian and $f$ will have zero integral at every time~$t$.
The Hamiltonian vector f\/ield is
$H_\ps=\partial_y\ps\partial_x-\partial_x\ps\partial_y$
and the Poisson bracket $L_{H_\ps}f=\{f,\ps\}$
is the Jacobian of $f$ and $\ps$.
If $\al=-\be dy$ and $a=-1$ we get the equation
for stream function $\ps$ and vorticity function $\om=\De\ps$:
\begin{gather}
\partial_t \om =-\{\om,\ps\}-\be\partial_x f,\nonumber\\
\partial_t f =-\{f,\ps\}+\be\partial_x\ps.\label{stratified}
\end{gather}

Let $\xi=g\frac{\rh-\rh_0}{\rh_0}$ be a buoyancy variable
measuring the deviation of a density $\rh$ from a background value~$\rh_0$,
with $g$ the gravity acceleration.
The background stratif\/ication $\rh_0$ is assumed to be exponential,
characterized by the constant Brunt--V\"ais\"al\"a
frequency $N=(-g\frac{d\log\rh_0}{dy})^{\frac12}$.
The equation for a {\bf stratif\/ied f\/luid in Boussinesq approximation}
\cite{Zeitlin2} is the geodesic equation~(\ref{stratified})
for $\be=N$ constant and $f=N^{-1}\xi$:
\begin{gather*}
\partial_t\om =-\{\om,\ps\}-\partial_x\xi,\\
\partial_t\xi =-\{\xi,\ps\}+N^2\partial_x \ps.
\end{gather*}

When the Brunt--V\"ais\"al\"a frequency $N$ is an integer and $\xi$
has zero integral (at time zero),
then the stratif\/ied f\/luid equation is a geodesic equation on a Lie group~\cite{V2}.


\section[$H^1$ metrics]{$\boldsymbol{H^1}$ metrics}

\subsection*{Camassa--Holm equation}

The {\bf Camassa--Holm equation} \cite{CH}
\begin{gather}\label{ch}
\partial_t(u-u'')=-3uu'+2u'u''+uu'''
\end{gather}
is the geodesic equation for the
right invariant metric on $\Diff(S^1)$ given by the $H^1$ scalar product
$\langle X,Y\rangle =\int_{S^1}(XY+X'Y')dx=\int_{S^1}X(1-\partial_x^2)Ydx$
\cite{Kouranbaeva}.

Indeed, one gets from
\begin{gather*}
\langle \ad(X)^\top Y,Z\rangle  =\langle Y,X'Z-XZ'\rangle
=\int_{S^1}(Y(X'Z-XZ')+Y'(X''Z-XZ''))dx\\
\phantom{\langle \ad(X)^\top Y,Z\rangle}{} =\int_{S^1}Z
(2YX'+Y'X-2Y''X'-Y'''X)dx
\end{gather*}
that $\ad(X)^\top Y=(1-\partial_x^2)^{-1}(2YX'+Y'X-2Y''X'-Y'''X)$.
Plugging
$\ad(X)^\top X=(1-\partial_x^2)^{-1}(3XX'-2X'X''-XX''')$
into Euler's equation (\ref{euler}) one obtains
the Camassa--Holm shallow water equation for $u:I\to C^\oo(S^1)$.

Since $m=A(u)=u-u''$,
the Hamiltonian form of the Camassa--Holm equation is
\begin{gather*}
\partial_tm=-um'-2u'm.
\end{gather*}

\begin{remark}
Considering the right invariant $H^1$ metric on the Bott--Virasoro group (\ref{bott}), an extended Camassa--Holm equation is obtained \cite{Misiolek1}
\begin{gather*}
\partial_t(u-u'')=-3uu'+2u'u''+uu'''-2au''',\qquad a\in\RR.
\end{gather*}
Indeed, the identity $\om(X,Y)=\langle k(X),Y\rangle$ for the Virasoro cocycle $\om(X,Y)=2\int_{S^1}X'''Ydx$ and the $H^1$ scalar product implies $k(X)=2(1-\partial_x^2)^{-1}X'''$. Now by
Corollary \ref{central} the geodesic equation
is the extended Camassa--Holm equation above.

The homogeneous manifold $\Diff(S^1)/S^1$ is a coadjoint orbit of the Bott--Virasoro group.
The {\bf Hunter--Saxton equation} describing weakly nonlinear unidirectional waves \cite{HS}
\begin{gather*}
\partial_tu''=-2u'u''-uu'''
\end{gather*}
is a geodesic equation on  $\Diff(S^1)/S^1$ with the right invariant metric def\/ined by the scalar product $\langle X,Y\rangle=\int_{S^1}X'Y'dx$ \cite{KM2}.
\end{remark}

\subsection*{Higher dimensional Camassa--Holm equation}

The {\bf higher dimensional Camassa--Holm equation} (also called EPDif\/f or
averaged template matching equation)
\cite{HMR,HM}
is the geodesic equation for the
right invariant $H^1$ metric
\begin{gather}\label{def}
\langle X,Y\rangle=\int_M(g(X,Y)+\al^2g(\nabla X,\nabla Y))\mu,
\end{gather}
on $\Diff(M)$ for compact $M$.
Because $\nabla^*\nabla=\De+\Ric$,
this scalar product can be rewritten with the help of the rough Laplacian
$\De_{R}=\De+\Ric$ as
$\langle X,Y\rangle=\int_M(g(X-\al^2\De_{R}X,Y)\mu$,
so the momentum density of the f\/luid $m=A(u)$ is $m=(1-\al^2\De_{R})u$.
It follows that the adjoint of $\ad(X)$ with respect to (\ref{def}) is conjugate
by $1-\al^2\De_R$ to the adjoint of $\ad(X)$ with respect to the $L^2$ metric (\ref{el2})
computed to be (\ref{brasov}).
Hence
$(1-\al^2\De_R)\ad(X)^\top Y=(\nabla X)^\top(Y-\al^2\De_RY)
+\nabla_X(Y-\al^2\De_RY)+(\div X)(Y-\al^2\De_RY)$.
We get as geodesic equation the higher dimensional Camassa--Holm equation
\begin{gather*}
\partial_t(1-\al^2\De_{R})u=-u\div u+\al^2(\div u)\De_{R} u-\nabla_uu
+\al^2\nabla_u(\De_{R} u)\\
\phantom{\partial_t(1-\al^2\De_{R})u=}{}
-(\nabla u)^\top u+\al^2(\nabla u)^\top \De_{R} u.
\end{gather*}
In Hamiltonian form this equation is
\begin{gather*}
\partial_tm=-\nabla_um-(\nabla u)^\top m-(\div u)m \qquad \mbox{for}\quad m=A(u)=u-\al^2\De_R u.
\end{gather*}

In particular for $M=S^1$ we get the Camassa--Holm equation (\ref{ch}).

When $M$ is a manifold with boundary and we put Neumann or mixed conditions
on the boundary, then the $H^1$ scalar product (\ref{def}) has to be replaced by
\begin{gather}\label{200}
\langle X,Y\rangle=\int_M(g(X,Y)+2\al^2g(\Def X,\Def Y))\mu,
\end{gather}
where $\Def X=\frac12(\nabla X+(\nabla X)^\top)$ denotes
the deformation $(1,1)$-tensor of $X$ \cite{GR2}.

\subsection*{Averaged Euler equation}

For a compact Riemannian manifold $M$ we consider the
right invariant metric on the group $\Diff_\mu(M)$
of volume preserving dif\/feomorphisms given by the $H^1$ scalar product
(\ref{def}) on vector f\/ields.
The geodesic equation is the (Lagrangian)
{\bf averaged Euler equation} \cite{MRS,Shkoller},
also called LAE-$\al$ equation:
\begin{gather}\label{lambda}
\partial_t(1-\al^2\De_R)u=-\nabla_u(1-\al^2\De_R)u
+\al^2(\nabla u)^\top(\De_R u)-\grad p.
\end{gather}
Indeed, from
$\langle\ad(X)^\top Y,Z\rangle=\int_Mg(Y-\al^2\De_RY,\nabla_ZX-\nabla_XZ)\mu
=\int_Mg((\nabla X)^\top (Y-\al^2\De_RY)+\nabla_X(Y-\al^2\De_RY),Z)\mu$
we obtain that
\begin{gather*}
(1-\al^2\De_R)(\ad(X)^\top Y)
=P((\nabla X)^\top Y+\nabla_X(1-\al^2\De_R)Y-\al^2(\nabla X)^\top(\De_RY))
\end{gather*}
and we use Euler's equation (\ref{euler}) to get (\ref{lambda}).
In Hamiltonian form this equation is
\begin{gather*}
\partial_tm=-\nabla_um-(\nabla u)^\top m-\grad p
\end{gather*}
for $m=A(u)=u-\al^2\De_R u$.

As in the higher dimensional Camassa--Holm equation,
when Neumann or mixed conditions on the boundary of $M$
are imposed, one has to consider the $H^1$ scalar product (\ref{200}).


\section{Systems of two evolutionary equations}

From \cite{Fuks} we know that a basis for $H^2(\mathfrak{X}(S^1),C^\oo(S^1))$ is represented by
$\si(X,Y)=X'Y-XY'$ and the Virasoro cocycle
$\om(X,Y)=\int_{S^1}(X'Y''-X''Y')dx\in\RR\subset C^\oo(S^1)$; a basis for
$H^2(\mathfrak{X}(S^1),\Om^1(S^1))$ is represented by the cocycles
\begin{gather*}
\si_1(X,Y)=XY''-X''Y,\qquad
\om_1(X,Y)=X'Y''-X''Y';
\end{gather*}
a basis for $H^2(\mathfrak{X}(S^1),\Om^2(S^1))$ is represented by the cocycles
\begin{gather*}
\si_2(X,Y)=X'''Y-XY''',\qquad
\om_2(X,Y)=X'''Y'-X'Y'''.
\end{gather*}
Only the cocycles $\om$, $\om_1$ and $\om_2$
(whose expressions involve only derivatives of $X$ and $Y$)
integrate to group cocycles \cite{OR}.

The Euler equations for the $L^2$ or $H^1$ scalar product on the corresponding
Abelian extensions provide systems of two equations,
generalizing Burgers (\ref{burgers}) or Camassa--Holm (\ref{ch}) equation.
We exemplify with the 2-cocycle $\si$ taking values in the module of functions on the circle.
The Euler equations for the $L^2$ scalar product on $C^\oo(S^1)\rtimes\X(S^1)$
and on $C^\oo(S^1)\rtimes_\si\X(S^1)$ are
\begin{gather*}
\partial_tu =-3uu'-ff',\\
\partial_tf =-uf'-u'f,
\end{gather*}
and
\begin{gather*}
\partial_tu =-3uu'+uf'+2u'f-ff',\\
\partial_tf =-uf'-u'f.
\end{gather*}
The Euler equations for the $H^1$ scalar product on $C^\oo(S^1)\rtimes\X(S^1)$
and on $C^\oo(S^1)\rtimes_\si\X(S^1)$ are
\begin{gather*}
\partial_t(u-u'') =-3uu'+2u'u''+uu'''-ff'+f'f''',\\
\partial_t(f-f'') =-uf'-u'f+uf'''+u'f'',
\end{gather*}
and
\begin{gather*}
\partial_t(u-u'') =-3uu'+2u'u''+uu'''-2u'f-uf'+2u'f''+uf'''-ff'+f'f''',\\
\partial_t(f-f'') =-uf'-u'f+uf'''+u'f''.
\end{gather*}

One can consider central extensions of semidirect products
of $\X(S^1)$ with modules of densities as in Remark~\ref{ov} \cite{OR}.
For instance the 1-cocycle $\al(X)=X''$ on $\X(S^1)$ with values in
$\Om^1(S^1)$, the module dual to $C^\oo(S^1)$,
gives the 2-cocycle $\om((f_1,X_1),(f_2,X_2))=\int_{S^1}(X_1''f_2-X_2''f_1)dx$
on the semidirect product $C^\oo(S^1)\rtimes\X(S^1)$.
The geodesic equation for the $L^2$ scalar product on the central extension
$(C^\oo(S^1)\rtimes\X(S^1))\times_\om \RR$ is
\begin{gather*}
\partial_tu =-3uu'-ff'-af'',\\
\partial_tf =-uf'-u'f-au'',\quad a\in\RR.
\end{gather*}


\section{Conclusions}
This survey article presents the formal deduction as geodesic equations on dif\/feomorphism groups with right invariant metrics of several PDE's of hydrodynamical type. Sometimes extensions of dif\/feomorphism groups by central or Abelian sugroups come into play and the corresponding Lie algebra 2-cocycles introduce additional terms to the geodesic equations.

These equations are Hamiltonian equations too, possessing rich geometric structures. Some of them are completely integrable. But presenting these results is beyond the scope of this article.


\subsection*{Acknowledgements}
This work was done with the f\/inancial support of Romanian Ministery of Education
and Research under the grant CNCSIS 95GR/2007.
I acknowledge the support from the ICTP Of\/f\/ice of External Activities for attending
the Seventh International Conference ``Symmetry in Nonlinear Mathematical Physics'' in Kyiv.

I am most grateful to Tudor Ratiu and Francois Gay-Balmaz for their preprints
and for very good suggestions, and to the referees for their very substantial and
constructive comments.

\pdfbookmark[1]{References}{ref}
\LastPageEnding

\end{document}